\documentclass[12pt]{amsart} 
\usepackage[mathscr]{eucal}
\usepackage{amsmath,amsfonts}

\parskip=\smallskipamount

\newtheorem{theorem}{Theorem}

\newtheorem{proposition}[theorem]{Proposition}

\newtheorem{lemma}[theorem]{Lemma}

\makeatletter
\@addtoreset{equation}{section}
\makeatother

\newcommand{\cA}{{\mathcal A}}
\newcommand{\cB}{{\mathcal B}}
\newcommand{\cC}{{\mathcal C}}
\newcommand{\cD}{{\mathcal D}}
\newcommand{\cE}{{\mathcal E}}
\newcommand{\cF}{{\mathcal F}}
\newcommand{\cG}{{\mathcal G}}
\newcommand{\cH}{{\mathcal H}}

\newcommand{\cL}{{\mathcal L}}

\newcommand{\cN}{{\mathcal N}}
\newcommand{\cP}{{\mathcal P}}
\newcommand{\cR}{{\mathcal R}}
\newcommand{\cS}{{\mathcal S}}

\newcommand{\cO}{{\mathcal O}}
\newcommand{\cI}{{\mathcal I}}

\newcommand{\CC}{{\mathbb C}}
\newcommand{\NN}{{\mathbb N}}

\newcommand{\RR}{{\mathbb R}}
\newcommand{\FF}{{\mathbb F}}
\newcommand{\TT}{{\mathbb T}}

\newdimen\expt
\def\boxit#1{\setbox0\hbox{$\displaystyle{#1}$}
      \hbox{\lower.4\expt
 \hbox{\lower3\expt\hbox{\lower\dp0
      \hbox{\vbox{\hrule height.4\expt
 \hbox{\vrule width.4\expt\hskip3\expt
      \vbox{\vskip3\expt\box0\vskip2\expt}%
 \hskip3\expt\vrule width.4\expt}\hrule height.4\expt}}}}}}

\begin{document}
\pagestyle{plain}

\bigskip

\title 
{Orthogonal polynomials  \\
in several variables. I} 
\author{T. Constantinescu} 

\address{Department of Mathematics \\
  University of Texas at Dallas \\
  Box 830688, Richardson, TX 75083-0688, U. S. A.}
\email{\tt tiberiu@utdallas.edu}

\maketitle

\begin{abstract}
In this paper we introduce and discuss some classes of orthogonal 
polynomials in several non-commuting variables. The emphasis 
is on a non-commutative version of the orthogonal
polynomials on the real line. We introduce recurrence equations
for these polynomials, 
Christoffel-Darboux formulas, and Jacobi type matrices.
\end{abstract}

\section{Introduction}

Orthogonal polynomials in several variables are known for long time,
see for instance \cite{EMOT}, but their theory is less developed
than in the one variable case. The commutative case is also 
studied more intensively
(see, for instance, \cite{BCR}, \cite{Xu}), while the studies
for the non-commutative case appear to be quite sparse
(see \cite{Ko}). 

Our goal is to introduce and study some classes of orthogonal 
polynomials in several non-commuting variables. In this paper we focus
on polynomials that are viewed as analogues of the orthogonal 
polynomials on the real line. The main topics are:
recurrence equations, Szeg\"o kernels, and Jacobi matrices. 

The paper is organized as follows. In Section~2 we introduce the main 
definitions and several examples. Especially, we briefly review
a non-commutative version of the Szeg\"o theory of orthogonal polynomials
on the unit circle that was sketched in \cite{CJ2}. Section~3 deals
with recurrence equations. In Section~4 we introduce 
a non-commutative Szeg\"o type kernel which is viewed as a sort
of reproducing kernel for the Siegel
upper half-space. Section~5 deals with Jacobi type matrices associated
to the recurrence equations introduced in Section~3. 

\section{Orthogonal polynomials}
In this section we introduce the main definitions and briefly discuss 
several examples. 
Let $\FF _M^+$ be the unital free semigroup 
on $M$ generators $g_1,\ldots ,g_M$ with lexicograhpic
order $\prec $. The empty word is the identity element
and the length of the word $\sigma $ is 
denoted by $|\sigma |$. The length of the empty word is $0$.

Let $\cO _N^0$ be the algebra of polynomials
in $2N$ non-commuting
indeterminates $Y_1$,$\ldots $,$Y_N$, $Y_{N+1}$,$\ldots $,$Y_{2N}$ 
with complex coefficients.
Each element $Q\in \cO _N^0$ can be uniquely written 
in the form $Q=\sum _{\sigma \in \FF _{2N}^+}c_{\sigma }Y_{\sigma }$,
with only finitely many $c_{\sigma }\ne 0$ and 
$Y_{\sigma }=Y_{i_1}\ldots Y_{i_k}$ for $\sigma =i_1\ldots i_k
\in \FF _{2N}^+$. An involution $\cI$ can 
be introduced on $\cO _N^0$ as follows:
$$\cI (Y_k)=Y_{N+k},\quad k=1,\ldots ,N,$$
$$\cI (Y_l)=Y_{l-N},\quad l=N+1,\ldots ,2N;$$
on monomials, 
$$\cI (Y_{i_1}\ldots Y_{i_k})=\cI (Y_{i_k})\ldots \cI (Y_{i_1}),$$
and finally, if $Q=\sum _{\sigma \in \FF _{2N}^+}c_{\sigma }Y_{\sigma }$,
then $\cI (Q)=
\sum _{\sigma \in \FF _{2N}^+}\overline{c}_{\sigma }
\cI (Y_{\sigma })$.
Thus, $\cO _N^0$ is a unital, associative, $*$-algebra over $\CC $.

Let $\cP _N^0$ denote the algebra of polynomials in $N$ non-commuting
indeterminates $Y_1,\ldots ,Y_N$ with complex coefficients.
Then $\cP _N^0$ is a subalgebra of $\cO _N^0$.
We say that $\cA \subset \cO _N^0$
is $\cI $-symmetric if $P\in \cA $ implies 
$c\cI (P)\in \cA $ for some $c\in \CC -\{0\}$. 
We construct an associative algebra
$\cO _N^0(\cA )$ as the quotient of $\cO _N^0$
by the two-sided ideal $\cE (\cA )$ generated by 
$\cA $. We notice that $\cO _N^0(\emptyset )=
\cO _N^0$.  We let $\pi =\pi _{\cA }:\cO _N^0\rightarrow 
\cO _N^0(\cA )$ be the quotient map and since
$\cA $ is $\cI $-symmetric,
\begin{equation}\label{invo}
\cI _{\cA }(\pi (P))=\pi (\cI (P))
\end{equation}
gives an involution on $\cO _N^0(\cA )$.
We will be interested in linear functionals $\phi $
on $\cO _N^0(\cA )$ with the property that 
$\phi (\cI _{\cA }(\pi (P))\pi (P))\geq 0$
for all $P\in \cP _N^0$. Without loss of generality we will
assume that $\phi $ is unital, $\phi (\pi (1))=1$. 
Such a functional will be called a 
{\em positive functional} on $\cO _N^0(\cA )$.
The proof of the following result is straightforward and 
can be ommited.

\begin{lemma}\label{1}
Let $\phi $ be a positive functional on $\cO _N^0(\cA )$. Then

$$\begin{array}{l}
1)\quad  \phi (\cI _{\cA }(\pi (P)))=
\overline{\phi (\pi (P))}\quad \mbox{for $P\in \cP _N^0$}.
\\
 \\
2)\quad  |\phi (\cI _{\cA }(\pi (P_1))\pi (P_2))|^2\leq 
\phi (\cI _{\cA }(\pi (P_1))\pi (P_1))
\phi (\cI _{\cA }(\pi (P_2))\pi (P_2))
\quad \mbox{for $P_1,P_2\in \cP _N^0$}.
\end{array}
$$
\end{lemma}

We now consider the GNS construction associated to
$\phi $. Thus, we define 
on $\pi (\cP _N^0)$,
\begin{equation}\label{tilde}
\langle \pi (P_1),\pi (P_2)\rangle _{\phi }=
\phi (\cI _{\cA }(\pi (P_2)\pi (P_1)),
\end{equation}
and factor out the subspace
$\cN _{\phi }=\{\pi (P)\mid P\in \cP _N^0,
\langle \pi (P),\pi (P) \rangle _{\phi }=0\}$.
Completing this quotient with respect
to the norm induced by \eqref{tilde} we obtain a Hilbert space 
$\cH _{\phi }$. From now on we will assume that 
$\phi $ 
is strictly positive, that is, 
$\phi (\cI _{\cA }(\pi(P))\pi (P))>0$
for all $P\in \cP _N^0-\cE (\cA)$, so that  $\cN _{\phi }=\{0\}$ and 
$\pi (\cP _N^0)$ can be viewed as a subspace of $\cH _{\phi }$.
Let $\cF =\{F_{\alpha }\}_{\alpha \in G}$ be the set 
of the distinct elements $\pi (Y_{\sigma })$, $\sigma \in 
\FF _N^+$. The {\em index set } $G$ is  
chosen as follows:
the equality on $\{\pi (Y_{\sigma })\mid \sigma \in \FF _N^+\}$
gives an equivalence relation and choose from each equivalence
class the element $\pi (Y_{\sigma })$ with the least $\sigma $
with respect to the lexicographic order. Then $G\subset \FF _N^+$
and $\emptyset \in G$.
Let $G_n=\{\alpha \in G\mid |\alpha |=n\}$, then $G_0=
\{\emptyset \}$ and $\{G_n\}_{n\geq 0}$ is a partition of $G$.

Since $\phi $ is strictly positive it follows that 
$\cF $ is a linearly independent family in 
$\cH _{\phi }$ and 
the 
Gram-Schmidt procedure gives a family $\{\varphi _{\alpha  }\}
_{\alpha \in G}$ of elements in 
$\pi (\cP ^0_N)\subset \cO _N^0(\cA )$ such that 

\begin{equation}\label{bond1}
\varphi _{\alpha  }=
\sum _{\beta \preceq \alpha }a_{\alpha ,\beta }F_{\beta },
\quad a_{\alpha ,\alpha }>0;
\end{equation}
\begin{equation}\label{bond2}
\langle \varphi _{\alpha }, \varphi _{\beta }\rangle _{\phi }=
\delta _{\alpha ,\beta },
\quad \alpha ,\beta \in G.
\end{equation}

\noindent
The elements $\varphi _{\alpha  }$, $\alpha \in G$, will be called
the {\em orthogonal polynomials} associated to $\phi $.
Typically, the theory of orthogonal polynomials deals
with the study of algebraic and asymptotic properties of
the orthogonal polynomials 
associated to strictly positive
functionals on $\cO _N^0(\cA )$.
An explicit formula for the orthogonal polynomials  
can be obtained in the same manner as in 
the classical (one variable) case. Define
\begin{equation}\label{stilde}
s_{\alpha ,\beta }=
\phi (\cI _{\cA }(F_{\alpha })F_{\beta })=
\langle F_{\beta },F_{\alpha }\rangle _{\phi }, 
\quad \alpha ,\beta \in G,
\end{equation}
and 
\begin{equation}\label{Deh}
D_{\alpha }=\det\left[s_{\alpha ',\beta '}\right]_
{\alpha ',\beta '\preceq \alpha }>0, \quad \alpha \in G.
\end{equation}  
We notice that $\phi $ is a positive functional on 
$\cO _N^0(\cA )$ if and only if $K_{\phi }(\alpha ,\beta )=
s_{\alpha ,\beta }$, $\alpha ,\beta \in G$, 
is a positive definite kernel on $G$.
From now on $\tau -1$ denotes the predecessor of $\tau $ with respect to the 
lexicographic order $\prec $ on $\FF _N^+$, while $\sigma +1$
denotes the successor of $\sigma $. The determinant involved in the 
next result is defined by the same formula as in the scalar case,
even though its entries are elements of $\pi (\cP _N^0)$.

\begin{theorem}\label{T2}
Let $\{\varphi _{\alpha }\}_{\alpha \in G}$
be the orthogonal polynomials associated to the 
strictly positive, unital functional $\phi $
on $\cO _N^0(\cA )$.
Then $\varphi _{\emptyset }=1$ and for $\emptyset \prec \alpha $,
\begin{equation}\label{crop1}
\varphi _{\alpha }=\frac{1}{\sqrt{D_{\alpha -1}D_{\alpha }}}
{\det \left[\begin{array}{c}
\left[s_{\alpha ',\beta '}\right]_
{\alpha '\prec \alpha ;
\beta '\preceq \alpha } \\
  \\
\begin{array}{ccc}
F_{\emptyset } & \ldots & F_{\alpha }
\end{array}
\end{array}
\right]}.
\end{equation}
\end{theorem}

\begin{proof}
The proof is similar to the classical one. Thus, we deduce 
from the orthogonality condition \eqref{bond2} that
$\langle \varphi _{\alpha  },F_{\beta '}\rangle _{\phi }=0$
for $\emptyset \preceq \beta '\prec \alpha $, which implies 
that $\sum _{\beta \preceq \alpha }a_{\alpha ,\beta }
s_{\beta ',\beta }=0$ for $\emptyset \preceq \beta '\prec \alpha $.
Since the coefficients of the linear system 
$$\left\{\begin{array}{rcl}
\sum _{\beta \preceq \alpha }a_{\alpha ,\beta }
s_{\beta ',\beta }&=&0, \quad \emptyset \preceq \beta '\prec \alpha ,\\
 & & \\
\sum _{\beta \preceq \alpha }a_{\alpha ,\beta }F_{\beta }&=
&\varphi _{\alpha },
\end{array}
\right.
$$
with unknowns $a_{\alpha ,\beta }$
are complex
numbers except for those of the last equation which are
in $\pi (\cP _N^0)$, Cramer's  rule still holds 
in a form that gives 
$$\varphi _{\alpha }=\frac{a_{\alpha ,\alpha }}{D_{\alpha -1}}
\det \left[\begin{array}{c}
\left[s_{\alpha ',\beta '}\right]_
{\alpha '\prec \alpha ;\beta '\preceq \alpha } \\
\\
\begin{array}{ccc}
F_{\emptyset } & \ldots & F_{\alpha }
\end{array}
\end{array}
\right].$$
Next we notice that 
$$\langle \det \left[\begin{array}{c}
\left[s_{\alpha ',\beta '}\right]_
{\alpha '\prec \alpha ;\beta '\preceq \alpha } \\
\\
\begin{array}{cccc}
F_{\emptyset } & \ldots & F_{\alpha }
\end{array}
\end{array}
\right],F_{\alpha }\rangle _{\phi }=D_{\alpha  },$$
and since $F_{\alpha }=\frac{1}{a_{\alpha ,\alpha }}\varphi _{\alpha  }+
\sum _{\beta \prec \alpha }c_{\beta }F_{\beta }$ for some complex
coefficients $c_{\beta }$, $\beta \prec \alpha $,
we deduce 
$$D_{\alpha }=\langle \frac{D_{\alpha -1}}{a_{\alpha ,\alpha }}
\varphi _{\alpha },\frac{1}{a_{\alpha ,\alpha }}\varphi _{\alpha }+
\sum _{\beta \prec \alpha }c_{\beta }F_{\beta }\rangle _{\phi }=
\frac{D_{\alpha -1}}
{a^2_{\alpha ,\alpha }},$$
so that 
$$\frac{1}{a^2_{\alpha ,\alpha }}=\frac{D_{\alpha }}{D_{\alpha -1}},$$
which gives \eqref{crop1}. 
\end{proof}

Usually, the representation \eqref{crop1}
is not very useful for the actual computation of the orthogonal
polynomials. Instead, recurrence equations are obtained
for each particular case of interest. We consider several such examples.

\bigskip
\noindent
{\em 2.1. Orthogonal polynomials in one variable.}
Assume $N=1$. Then a linear functional $\phi $
on $\cO _1^0$ is positive if and only if 
$K_{\phi }(n,m)=\phi (\cI (Y_1^n)Y_1^m)$, $n,m\in \NN $, 
is a positive definite kernel on $\NN $.
One simple example can be obtained by taking
$\cA _{\TT }=\{1-\cI (Y_1)Y_1\}.$
In this case, 
$\phi $ is a positive functional on 
$\cO _1^0(\cA _{\TT})$ if and only if the kernel $K_{\phi }$
is positive definite and satisfy the Toeplitz condition, that is
$$K_{\phi }(n+k,m+k)=K_{\phi }(n,m), \quad m,n,k\in \NN .$$

This shows that the orthogonal polynomials on 
$\cO _1^0(\cA _{\TT})$
are the orthogonal polynomials on the unit circle (\cite{Sz}).
The index set is $\NN $ in this case, and define
$\gamma _n=-a_{n,n}\varphi _n(0)$ for $n\geq 1$. Then $|\gamma _n|<1$
and set $d_n=(1-|\gamma _n|^2)^{1/2}$. The orthogonal 
polynomials on the unit circle satisfy the following basic 
recurrence equations (see \cite{Sz}):

\begin{equation}\label{recur1}
\varphi _{n+1}(z)=\frac{1}{d_{n+1}}(z\varphi _n(z)-
\gamma _{n+1}\varphi _n^{\sharp }(z)),
\end{equation}
where $\varphi _0^{\sharp }(z)=1$ and for $n>0$, 

\begin{equation}\label{recur2}
\varphi _{n+1}^{\sharp }(z)=\frac{1}{d_{n+1}}
(-\overline{\gamma }_{n+1}z\varphi _n(z)+
\varphi _n^{\sharp }(z)).
\end{equation}

Another example is given by 
$\cA _{\RR }=\{Y_1-\cI (Y_1)\}$.
The index set is still $\NN $, 
but this time, $\phi $ is a positive functional on  
$\cO _1^0(\cA _{\RR})$
if and only if  
the kernel $K_{\phi }$
is positive definite and has the Hankel property, that is
$$K_{\phi }(n,m+k)=K_{\phi }(n+k,m), \quad m,n,k\in \NN .$$

This shows that the orthogonal polynomials on 
$\cO _1^0(\cA _{\RR})$ are the orthogonal polynomials on the 
real line. In this case one obtains a basic 
three-terms recurrence equation (see \cite{Sz}):

\begin{equation}\label{iarrecur}
x\varphi _n(x)=b_{n+1}\varphi _{n+1}(x)+
a_{n}\varphi _n(x)+
b_{n}\varphi _{n-1}(x),
\end{equation}

\noindent
with initial conditions $\varphi _0(x)=1$ and $\varphi _{-1}(x)=0$.

\bigskip
\noindent
{\em 2.2. Examples in several variables.}
In this case, there is a large number of interesting examples.
Here we mention just one. More examples will be considered in the 
next sections. 
Let 
$$\cA _1=\{1-\cI (Y_k)Y_k \mid k=1,\ldots N\}\cup
\{\cI (Y_k)Y_l\mid k,l=1,\ldots N, k\ne l\}.
$$
Then $\phi $ is a positive functional on 
$\cO _N^0(\cA _1)$ if and only if $K_{\phi }$
is a positive definite kernel that obeys the rules:

\begin{equation}\label{sta1}
K_{\tilde \phi }(\tau \sigma ,\tau \sigma ')=
K_{\tilde \phi }(\sigma ,\sigma '), \quad 
\tau ,\sigma ,\sigma '\in \FF _N^+,
\end{equation}
\begin{equation}\label{sta2}
K_{\tilde \phi }(\sigma ,\tau )=0 
\quad \mbox{ if there is no $\alpha \in \FF _N^+$
such that $\sigma =\alpha \tau $ or $\tau =\alpha \sigma $}.
\end{equation}

\noindent
Such type of kernels appeared in the study of some classes of 
stochastic processes indexed by nodes on a tree (see, 
for instance, \cite{CWB}). In this case the 
index set is $\FF _N^+$. We also need to introduce

\begin{equation}\label{1D}
D_{1,\sigma }=\det\left[K_{\sigma ',\tau '}\right]_
{\emptyset \prec \sigma ',\tau '\preceq \sigma }>0.
\end{equation} 

\noindent
Define $\gamma _{\sigma }=-\sqrt{\frac{D_{\sigma }}
{D_{1,\sigma }}}a_{\sigma ,\emptyset }$ and notice that
$|\gamma _{\sigma }|<1$. Then we can define 
$d_{\sigma }=(1-|\gamma _{\sigma }|^2)^{1/2}$ and 
the orhogonal 
polynomials associated to 
$\cO _N^0(\cA _1)$ obey the recurrence equations 
(see \cite{CJ2}):

\begin{equation}\label{szego}
\varphi _{k\sigma }=\frac{1}{d_{k\sigma }}
(Y_k\varphi _{\sigma }-\gamma _{k\sigma } 
\varphi ^{\sharp}_{k\sigma -1}),\quad k=1,\ldots ,N,\,\,
\sigma \in \FF _N^+,
\end{equation}
where $\varphi ^{\sharp}_{\emptyset }=1$ and for $k\in \{1,\ldots ,N\}$,
\,\,$\sigma \in \FF _N^+$,
\begin{equation}\label{sarp}
\varphi ^{\sharp}_{k\sigma }=\frac{1}{d_{k\sigma }}
(-\overline{\gamma }_{k\sigma }Y_k\varphi _{\sigma }+ 
\varphi ^{\sharp}_{k\sigma -1}).
\end{equation}

\section{Recurrence equations}
In this section we consider some algebraic properties 
of the orthogonal polynomials of 
$\cO _N^0(\cA _2)$, where 
$\cA _2=\{Y_k-\cI (Y_k)\mid k=1,\ldots ,N\}$.
If we take $\cA '_2=\cA _2\cup \{Y_kY_l-Y_lY_k
\mid k,l=1,\ldots ,N\}$, then $\pi _{\cA '_2}(\cP _N^0)$
is isomorphic to the symmetric algebra of $\CC ^N$
and the orthogonal polynomials correspond to the ortogonal polynomials
of several commuting variables, see \cite{Xu}. Since   
$\pi _{\cA '_2}(\cP _N^0)$ is a quotient of 
$\pi _{\cA _2}(\cP _N^0)=\cP _N^0$, 
we expect that results for $\cP _N^0$
would give corresponding results for 
$\pi _{\cA '_2}(\cP _N^0)$ by standard (symmetrization) techniques, 
see \cite{Pa}.
Going in the opposite direction, we expect to deduce generalizations
of results on $\pi _{\cA '_2}(\cP _N^0)$ to similar results for
$\cP _N^0$. We illustrate this remark by obtaining a three-term
recurrence relation for the orthogonal polynomials in 
$\cO _N^0(\cA _2)$. 

Let $\phi $ be a strictly positive functional on 
$\cO _N^0(\cA _2)$. Since $\pi _{\cA _2}(\cP _N^0)=\cP _N^0$,
the index set $G$ is $\FF _N^+$ and $G_n$ is the set of words
in $\FF _N^+$ of length $n$.
Let $\{\varphi _{\sigma }\}_{\sigma \in \FF _N^+}$
be the orthogonal polynomials associated to $\phi $.
The matrix-vector notation in \cite{Kow}
is easily adapted to $\cP _N^0$, by setting  
$\Phi _n=\left[\varphi _{\sigma }\right]_{|\sigma |=n},$
$n\geq 0$.
We can show that \eqref{iarrecur} extends to 
$\cO _N^0(\cA _2)$.

\begin{theorem}\label{T3}
The orthogonal polynomials on $\cO _N^0(\cA _2)$
obey the following recurrence relations: for $k=1,\ldots ,N$,
\begin{equation}\label{treiunu}
Y_k\Phi _0=
\Phi _{1}B_{0,k}+
\Phi _0A_{0,k},
\end{equation}
and for 
$k=1,\ldots ,N$ and $n\geq 1$,

\begin{equation}\label{rere}
Y_k\Phi _n=
\Phi _{n+1}B_{n,k}+
\Phi _{n}A_{n,k}+
\Phi _{n-1}B^*_{n-1,k}.
\end{equation}
\end{theorem}
\begin{proof}
The proof is most like in the classical, one-dimensional case.
Thus, we can write
$$Y_k\varphi _{\sigma }=\sum _{\tau \preceq k\sigma }
c^{k\sigma }_{\tau }\varphi _{\tau },\quad 
k\in \{1,\ldots ,N\},\,\,\sigma \in \FF _N^+,$$
where the coefficients $c^{k\sigma }_{\tau }$ are calculated by 
the formula $c^{k\sigma }_{\tau }=\langle Y_k
\varphi _{\sigma },\varphi _{\tau }\rangle _{\phi }$.
We notice that for any $P,Q\in \cP _N^0$,
\begin{equation}\label{simetrie}
\begin{array}{rcl}
\langle Y_kP,Q\rangle _{\phi }&=&\phi (\cI (Q)Y_kP) \\
 & & \\
&=&\phi (\cI (Q)\cI (Y_k)P) \\
 & & \\
&=&\langle P,Y_kQ\rangle _{\phi }.
\end{array}
\end{equation}

\noindent
[More generally, we have 
$s_{\alpha \sigma ,\tau }=s_{\sigma ,I (\alpha )\tau}$ for 
$\alpha ,\sigma ,\tau \in \FF _N^+$,
where $I $ denotes the involution on  
$\FF _N^+$ given by $I(i_1\ldots i_k)=i_k\ldots i_1$].
In particular, for $|\tau |\leq |\sigma |-2$,
$$c^{k\sigma }_{\tau }=\langle
\varphi _{\sigma },Y_k\varphi _{\tau }\rangle _{\phi }=0,$$
while for the remaining values of $\tau \preceq k\sigma $,
\begin{equation}\label{bine}
c^{k\sigma }_{\tau }=
\langle
\varphi _{\sigma },Y_k\varphi _{\tau }\rangle _{\phi }=
\overline{\langle Y_k\varphi _{\tau },
\varphi _{\sigma }\rangle _{\phi }}=
\overline{c}^{k\tau }_{\sigma }.
\end{equation}
We deduce that 
for $k=1,\ldots ,N$,
\begin{equation}\label{cere}
Y_k\Phi _0=
\Phi _1B_{0,k}+
\Phi _0A_{0,k},
\end{equation} 
while for 
$k=1,\ldots ,N$ and $n\geq 1$,

\begin{equation}\label{pere}
Y_k\Phi _{n}=
\Phi _{n+1}B_{n,k}+
\Phi _{n}A_{n,k}+
\Phi _{n-1}C_{n,k}.
\end{equation}

Let $X=\left[x_{ij}\right]$ be a given matrix.
We will use the following 
notation: first, 
$\cI (X)=\left[\cI (x_{ji})\right]$
and then $\phi (X)=\left[\phi (x_{ij})\right]$.
We deduce from \eqref{cere}, \eqref{pere}, and \eqref{bond2}
that 
\begin{equation}\label{aaa}
B^*_{n,k}=
\phi (\cI (\Phi _n)Y_k\Phi _{n+1}),
\end{equation}
\begin{equation}\label{bbb}
A^*_{n,k}=
\phi (\cI (\Phi _n)
Y_k\Phi _n),
\end{equation} 
and
\begin{equation}\label{ccc}
C^*_{n,k}=
\phi (\cI (\Phi _n)
Y_k\Phi _{n-1})=\phi (\cI (\Phi _{n-1})
Y_k\Phi _{n})^*=B_{n-1,k}.
\end{equation}
\end{proof}

We notice that \eqref{simetrie} implies that 
$A^*_{n,k}=A_{n,k}$ and if we define 
$B_n=\left[\begin{array}{ccc}
B_{n,1} & \ldots & B_{n,N}
\end{array}
\right]$, $n\geq 0$, then $B_n$ is an
$N^{n+1}\times N^{n+1}$ upper triangular invertible matrix.
Under these conditions we can prove a converse of Theorem~\ref{T3}.
This appears as a Favard type result and gives a construction of
strictly positive functionals on $\cO _N^0(\cA _2)$.

\begin{theorem}\label{T4}
Let $\varphi _{\sigma }=\sum _{\tau \preceq \sigma }
a_{\sigma ,\tau}Y_{\tau }$, $\sigma \in \FF _N^+$,
be elements in $\cP _N^0$ such that  $\varphi _{\emptyset }=1$ and 
$a_{\sigma ,\sigma }>0$.
Assume that there exist families 
$\{A_{n,k}\mid n\geq 0, \,\,k=1,\ldots ,N\}$, 
$\{B_{n,k}\mid n\geq 0, \,\,k=1,\ldots ,N\}$
of matrices such that $A^*_{n,k}=A_{n,k}$ for all 
$n\geq 0$ and $k=1,\ldots ,N$, 
$B_n=\left[\begin{array}{ccc}
B_{n,1} & \ldots & B_{n,N}
\end{array}
\right]$ is  an upper triangular invertible matrix for each $n\geq 0$,
for $k=1,\ldots ,N$,
\begin{equation}\label{patruunu}
Y_k\varphi _{\emptyset }=
\left[\varphi _{\sigma }\right]_{|\sigma |=1}B_{0,k}+
\varphi _{\emptyset }A_{0,k},
\end{equation}
and for 
$k=1,\ldots ,N$, $n\geq 1$,

\begin{equation}\label{patrudoi}
Y_k\left[\varphi _{\sigma }\right]_{|\sigma |=n}=
\left[\varphi _{\sigma }\right]^T_{|\sigma |=n+1}B_{n,k}+
\left[\varphi _{\sigma }\right]^T_{|\sigma |=n}A_{n,k}+
\left[\varphi _{\sigma }\right]^T_{|\sigma |=n-1}B^*_{n-1,k}.
\end{equation}

\noindent
Then there exists a strictly positive functional $\phi $
on $\cO _N^0(\cA _2)$ such that 
$\{\varphi _{\sigma }\}_{\sigma \in \FF _N^+}$
is the family of orthogonal polynomials associated to $\phi $.
\end{theorem}
\begin{proof}
Since $B_n$, $n\geq 0$, are invertible matrices, it follows
that $\{\varphi _{\sigma }\}_{\sigma \in \FF _N^+}$
is a linearly independent family in $\cP _N^0$
and formula \eqref{bond2}
suggests to define:

\begin{equation}\label{definitie}
\phi (1)=1\quad \mbox{and}\quad \phi (\varphi _{\sigma })=0
\quad \mbox{for}\quad \sigma \in \FF _N^+-\{\emptyset \}.
\end{equation}

\noindent
These relations uniquely determine a linear functional on 
$\cO _N^0(\cA _2)=\cP _N^0$. We will use again the 
matrix-vector notation,
$\Phi _n=\left[\varphi _{\sigma }\right]_{|\sigma |=n}$, $n\geq 0.$
Also, set $\Phi _{-1}=0$ and $B_{-1,k}=0$. Then \eqref{patruunu}
can be included in \eqref{patrudoi} for $n=0$.

We prove that 
\begin{equation}\label{gata}
\phi (\cI (\Phi _n)\Phi _m)=
\left\{
\begin{array}{cc}
0, & n\ne m;\\
I, & n=m,
\end{array}
\right.
\end{equation}
where $0$ and $I$ denote the zero, respectively the identity
matrix of a suitable dimension.
This will imply that $\phi $ is a strictly positive functional on 
$\cP _N^0$
and that the orthogonal polynomials 
associated to $\phi $ are precisely 
$\varphi _{\sigma }$, $\sigma \in \FF _N^+$.

We prove \eqref{gata} by induction on $m$. For $m=0$, 
$$\begin{array}{rcl}
\phi (\cI (\Phi _n)\Phi _0)&=&\Phi (\cI (\Phi _0)\Phi _n)^* \\
 & & \\
 &=&\left[\phi (\varphi _{\sigma })\right]^*_{|\sigma |=n},
\end{array}
$$
so that, by \eqref{definitie}, $\phi (\cI (\Phi _0)\Phi _0)=1$
and $\phi (\cI (\Phi _n)\Phi _0)=0$ for $n>0$.
Assume \eqref{gata} holds 
for $n\geq 0$ and $k\leq m$. Then $\phi (\cI (\Phi _n)\Phi _{m+1})=0$ 
for $n\leq m$ by the induction hypothesis. We deduce from 
\eqref{patrudoi}
that
$$\left[Y_1\Phi _l\ldots Y_N\Phi _l\right]=
\Phi _{l+1}B_l+\Phi _lA_l+\Phi _{l-1}C_l,$$
where $A_l=\left[A_{l,1}\ldots A_{l,N}\right]$
and $C_l=\left[B^*_{l-1,1}\ldots B^*_{l-1,N}\right]$, hence
$$\begin{array}{rcl}
\Phi _{l+1}&=&\left[Y_1\Phi _l\ldots Y_N\Phi _l\right]B^{-1}_l
-\Phi _lA_lB^{-1}_l-\Phi _{l-1}C_lB^{-1}_l \\
 & & \\
 &=&\sum _{j=1}^NY_j\Phi _lD_{l,j}-\Phi _lE_l
-\Phi _{l-1}F_l,
\end{array}
$$
where $B^{-1}_l=\left[D_{l,1}\ldots D_{l,N}\right]^T$,
$T$ denoting the matrix transpose,
$E_l=A_lB^{-1}_l$, and $F_l=C_lB^{-1}_l$.
Using the induction hypothesis, the previous formula for 
$l=m$, and \eqref{patrudoi}, we deduce
$$\begin{array}{rcl}
\phi (\cI (\Phi _{m+1})\Phi _{m+1})&=&
\sum _{j=1}^ND^*_{m,j}
\phi (\cI (\Phi _m)Y_j\Phi _{m+1}) \\
 & & \\ 
 & &-E^*_m\phi (\cI (\Phi _{m})\Phi _{m+1}) 
-F^*_m\phi (\cI (\Phi _{m-1})\Phi _{m+1}) \\
 & & \\
 &=&\sum _{j=1}^ND^*_{m,j}
\phi (\cI (Y_j\Phi _{m})\Phi _{m+1}) \\
 & & \\
 &=&\sum _{j=1}^ND^*_{m,j}
B^*_{m,j}\phi (\cI (\Phi _{m+1})\Phi _{m+1}) \\ 
 & & \\
 & &+\sum _{j=1}^ND^*_{m,j}B^*_{m,j}\phi (\cI (\Phi _{m})\Phi _{m+1}) \\
 & & \\
 & &+\sum _{j=1}^ND^*_{m,j}A_{m-1,j}\phi (\cI (\Phi _{m-1})\Phi _{m+1}) \\
 & & \\
 &=&\sum _{j=1}^ND^*_{m,j}B^*_{m,j}=I.
\end{array}
$$
Similar calculations show that $\phi (\cI (\Phi _{l})\Phi _{m+1})=0$ 
for $l>m+1$ and the proof is concluded.
\end{proof}

\section{Christoffel-Darboux formula}
In the classical theory orthogonal polynomials 
are evaluated at points in some suitable domains. For instance,
for orthogonal polynomials on the unit circle the domain is the 
open unit disk and the polynomials are closely related to the 
theory of analytic functions on that domain. A key to this connection
can be considered to be the Szeg\"o kernel $K_S(z,w)=\frac{1}{1-z
\overline{w}}$. A similar relation is established between orthogonal
polynomials on the real line and function theory on the 
upper half plane. In this section we suggest an extension of
these connections to several non-commuting variables. The analogue 
of the unit disk was already delt with in \cite{CJ1}, \cite{CJ2},
and we begin by briefly reviewing that construction.

 Let $\cE $ be 
an infinite-dimensional Hilbert space and
by $\cL (\cE )$ we denote the set of bounded linear operators
on $\cE $. The 
$N$-dimensional unit ball of $\cE $
is defined by
$$\cB _N(\cE )=\{Z=\left(\begin{array}{ccc}
Z_1 & \ldots & Z_N
\end{array}
\right)\mid (Z\mid Z)<I_{\cE}\},
$$
where $I_{\cE }$ is the identity on $\cE $ and 
for $Z=\left(\begin{array}{ccc}
Z_1 & \ldots & Z_N
\end{array}
\right)$ and 
$Z'=\left(\begin{array}{ccc}
Z'_1 & \ldots & Z'_N
\end{array}
\right)$
in $\cL (\cE )^N$,
\begin{equation}\label{inner}
(Z\mid Z')=\sum _{k=1}^NZ_k(Z')^*_k.
\end{equation}

A family of Hilbert spaces is associated to the Hilbert space $\cE $ 
as follows:
$\cE _0=\cE $ and for $k\geq 1$, 
$$\cE _k=\underbrace{\cE _{k-1}\oplus \ldots \oplus 
\cE _{k-1}}_{N\,\,terms}=\cE _{k-1}^{\oplus N}.$$
For $\cE =\CC $ we deduce $\CC _k=(\CC ^N)^{\otimes k}$, 
the $k$-fold tensor product of $\CC ^N$ with itself, therefore
$\oplus _{k\geq 0}\CC _k$ is the Fock space $\cF (\CC ^N)$
associated to $\CC ^N$. We also deduce that 
$\oplus _{k\geq 0}\cE _k$ is isomorphic to 
$\cF (\CC ^N)\otimes \cE $.

For $Z\in \cB _N(\cE )$
we define 
$
E(Z)=\left[Z_{\sigma }\right]_{|\sigma |=0}^{\infty }
$
and notice that $E(Z)$ is a bounded linear operator 
from $\oplus _{k\geq 0}\cE _k$ into $\cE$.
The {\em Szeg\"o kernel} for $\cB _N(\cE )$
is defined by the formula 
$$K_{\cB}(Z,Z')=E(Z)E(Z')^*,\quad Z,Z'\in \cB _N(\cE ).$$
The following result describes some of the basic properties of $K_{\cB}$.

\begin{lemma}\label{duminica}
$(a)$  $K_{\cB }$ is a positive definite kernel on  
$\cB _N(\cE )$.

\smallskip
\noindent
$(b)$ The set $\{E(Z)^*\cE \mid Z\in \cB _N(\cE )\}$
is total in $\oplus _{k\geq 0}\cE _k$.

\smallskip
\noindent
$(c)$ For any $T\in \cL (\cE )$ and $Z,Z'\in \cB _N(\cE )$, 
$$E(Z)(T-\sum _{k=1}^NZ_kT(Z')^*_k)^{\oplus \infty }E(Z')^*=T,$$
where $(T-\sum _{k=1}^NZ_kT(Z')^*_k)^{\oplus \infty }$
is the diagonal operator in $\cL (\oplus _{k\geq 0}\cE _k)$
with diagonal entry $T-\sum _{k=1}^NZ_kT(Z')^*_k$.
\end{lemma}
\begin{proof} The proof of this result can be found in 
\cite{CJ2}.   Thus, $(a)$ and $(c)$ are 
quite straightforward. The most interesting is $(b)$.
Its proof depends on the assumption the 
$\cE $ is infinite dimensional. In fact, for 
$\cE =\CC$, the result is not true, since the 
set $\{E(Z)^*\CC \mid Z\in \cB _N(\CC )\}$
is total in the symmetric Fock space, see \cite{Ar}.
For this reason and sake of completeness, we sketch the proof 
of $(b)$ here.
Let $f=\{f_{\sigma }\}_{\sigma \in \FF ^+_N}$
be an element of $\oplus _{k\geq 0}\cE _k$ orthogonal
to the linear span of $\{E(Z)^*\cE \mid Z\in \cB _N(\cE )\}$.
Taking $Z=0$, we deduce that $f_{\emptyset }=0$. Next, we claim
that for each $\sigma \in \FF ^+_N-\{\emptyset \}$ there exist
$$Z_l=(Z_1^l,\ldots ,Z_N^l)\in \cB _N(\cE ), \quad 
l=1,\ldots ,2|\sigma |,$$
such that 
$$\mbox{range}\left[
\begin{array}{ccc}
Z_{\sigma }^{*1} & \ldots & Z_{\sigma }^{*2|\sigma |}
\end{array}\right]=\cE,$$
and 
$$Z_{\tau }^l=0\quad \mbox{for all}\quad 
\tau \ne \sigma , \quad |\tau |\geq |\sigma |, \quad l=1,\ldots ,2|\sigma |.$$

Once this claim is proved, a simple inductive argument gives
$f=0$, so 
$\{E(Z)^*\cE \mid Z\in \cB _N(\cE )\}$
is total in $\oplus _{k\geq 0}\cE _k$.
In order to prove the claim we need the following construction.

Let $\{e^n_{ij}\}_{i,j=1}^n$ be the matrix units of the algebra
$M_n$ of $n\times n$ matrices. Each 
$e^n_{ij}$ is an $n\times n$ matrix consisting of $1$
in the $(i,j)th$ entry and zeros elsewhere.
For a Hilbert space $\cE _1$ we define $E^n_{ij}=e^n_{ij}\otimes 
I_{\cE _1}$ and we notice that 
\begin{equation}\label{units}
E^n_{ij}E^n_{kl}=\delta _{jk}E^n_{il},\quad E^{*n}_{ji}=E^n_{ij}.
\end{equation}

Let $\sigma =i_1\ldots i_k$ so that $\cE =\cE _1^{\oplus 2|\sigma |}$
for some Hilbert space $\cE _1$ (here we essentially use the assumption 
that $\cE $ is of infinite dimension).
Also, for $s=1,\ldots ,N$, we define 
$$J_s=\{l\in \{1,\ldots ,k\}\mid i_{k+1-l}=s\}$$
and 
$$
Z^{*p}_{s }=\frac{1}{\sqrt{2}}
\sum _{r\in J_s}E^{2|\sigma |}_{r+p-1,r+p},\quad 
s=1,\ldots ,N, \quad p=1,\ldots ,|\sigma |.$$
We can show that for each $p\in \{1,\ldots ,|\sigma \}$,
\begin{equation}\label{kiwi}
Z^{*p}_{\sigma }=\frac{1}{{\sqrt{2^k}}}E^{2|\sigma |}_{p,k+p},
\end{equation}
\begin{equation}\label{mango}
Z_{\tau }^p=0\quad \mbox{for} \quad 
\tau \ne \sigma ,\quad |\tau |\geq |\sigma |.
\end{equation}
Using
\eqref{units}, we deduce 
$$\begin{array}{rcl}
\sum _{s=1}^NZ_s^pZ_s^{*p}&=&
\frac{1}{2}\sum _{s=1}^N\sum _{r\in J_s}E^{2|\sigma |}_{r+p,r+p-1}
E^{2|\sigma |}_{r+p-1,r+p} \\
 & & \\
 &=&\frac{1}{2}\sum _{s=1}^N\sum _{r\in J_s}E^{2|\sigma |}_{r+p,r+p} \\
 & & \\
 &=&\frac{1}{2}\sum _{r=1}^kE^{2|\sigma |}_{r+p,r+p}<I,
\end{array}
$$
hence $Z^p\in \cB _N(\cE )$ for each $p=1,\ldots ,|\sigma |$.
For each word $\tau =j_1\ldots j_k\in \FF _N^+-\{\emptyset \}$
we deduce by induction that 
\begin{equation}\label{apple}
Z^{*p}_{j_k}\ldots 
Z^{*p}_{j_1}=\frac{1}{\sqrt{2^k}}
\sum _{r\in A_{\tau }}E^{2|\sigma |}_{r+p-1,r+p+k-1},
\end{equation}
where 
$A_{\tau }=\cap _{p=0}^{k-1}(J_{j_{k-p}}-p)\subset \{1,\ldots ,N\}$
and $J_{j_{k-p}}-p=\{l-p\mid l\in J_{i_{k-p}}\}$.

We show that $A_{\sigma }=\{1\}$ and $A_{\tau }=\emptyset $
for $\tau \ne \sigma $. Let $q\in A_{\tau }$. Therefore, for any 
$p\in \{0,\ldots ,k-1\}$ we must have $q+p\in J_{j_{k-p}}$
or $i_{k+1-q-p}=j_{k-p}$. For $p=k-1$ we deduce 
$j_1=i_{2-q}$ and since $2-q\geq 1$, it follows that $q\leq 1$. 
Also $q\geq 1$, therefore the only element that 
can be in $A_{\tau }$ is $q=1$, in which case we must have
$\tau =\sigma $. Since $l\in J_{i_{k+1-l}}$
for each $l=1, \ldots ,k-1$, hence $A_{\sigma }=\{1\}$
and $A_{\tau }=\emptyset $ for $\tau \ne \sigma $.
Formula \eqref{apple} implies \eqref{kiwi}. In a similar manner
we can construct a family $Z^p$, $p=|\sigma |+1, \ldots ,2|\sigma |,$
such that 
$$Z^{*p}_{\sigma }=\frac{1}{\sqrt{2^k}}E^{2|\sigma |}_{p+k,p},$$
and 
$$
Z_{\tau }^p=0\quad \mbox{for} \quad 
\tau \ne \sigma ,\quad |\tau |\geq |\sigma |.
$$
Thus, for $s=1,\ldots ,N$, we define 
$$K_s=\{l\in \{1,\ldots ,k\}\mid i_k=s\}$$
and 
$$Z^{*p}_{s}=\frac{1}{\sqrt{2}}
\sum _{r\in K_s}E^{2|\sigma |}_{r+p-k,r+p-k-1},\quad 
s=1,\ldots ,N, \quad p=|\sigma |+1,\ldots ,2|\sigma |.$$

Now, 
$$\left[
\begin{array}{ccc}
Z_{\sigma }^{*1} & \ldots & 
Z_{\sigma }^{*2|\sigma |}
\end{array}
\right]=
\frac{1}{\sqrt{2^k}}\left[
\begin{array}{cccccc}
E^{2|\sigma |}_{1,k+1} & \ldots & 
E^{2|\sigma |}_{k,2k} & E^{2|\sigma |}_{k+1,1}
& \ldots & E^{2|\sigma |}_{2k,k}
\end{array}
\right],
$$
whose range is $\cE $. 
This concludes the proof. 
\end{proof}

Define 
the vector space
$$\cR =\{r_f:\cB _N(\cE )\rightarrow \cE \mid 
r_f(Z)=E(Z)f, \,\, f\in \oplus _{k\geq 0}\cE _k\},$$
and consider the map $U:\oplus _{k\geq 0}\cE _k\rightarrow \cR$
defined by $Uf=r_f$. This map is linear and bijective (by Lemma 
~\ref{duminica} $(b)$), so we can define
on $\cR$ an inner product by the formula
$$\langle r_f,r_h\rangle _{\cR }=
\langle f,h\rangle _{\oplus _{k\geq 0}\cE _k}.$$
$\cR$ becomes a Hilbert space and $U$ is a unitary operator from 
$\oplus  _{k\geq 0}\cE _k$ onto $\cR $. The space $\cR $ has the 
{\em reproducing property}
\begin{equation}\label{repro}
\langle r_f(Z),e\rangle _{\cE }=
\langle r_f,r_{E(Z)^*e}\rangle _{\cR },
\quad f\in \oplus _{k\geq 0}\cE _k, \,\,e\in \cE .
\end{equation}
We now let $P=\sum _{\sigma \in \FF _N^+}c_{\sigma }Y_{\sigma }
\in \cP ^0_N$ take values on $\cB _N(\cE )$
by the the formula:
$$
P(Z)=\sum _{\sigma \in \FF _N^+}c_{\sigma }Z_{\sigma },
\quad 
Z\in \cB _N(\cE ).
$$
Note that each $f_e=\left[c_{\sigma }\right]_{\sigma \in \FF _N^+}
\oplus e$, $e\in \cE $, belongs to $\oplus _{k\geq 0}\cE _k$,
so that we deduce from \eqref{repro} that
$$\langle P(Z)e,e'\rangle _{\cE }=
\langle r_{f_e},r_{E(Z)^*e'}\rangle _{\cR },\quad e,e'\in \cE .$$

All these considerations suggest to introduce the following
set as a convenient domain for the 
evaluation of the orthogonal polynomials on 
$\cO _N^0(\cA _2)$.
Define 
$$\cG _N(\cE )=
\{\left(W_1\ldots W_N\right)\in \cL (\cE )^N\mid W_1W^*_1+\ldots 
+W_{N-1}W^*_{N-1}<\frac{1}{2i}(W_N-W^*_N)\},
$$
the {\em Siegel upper half-space} of $\cE $ (see \cite{Fo}, 
\cite{Ra}). There is a linear fractional map ({\em Cayley transform})
connecting $\cG _N(\cE )$ with $\cB _N(\cE )$, defined by the formula:
$$\cC (Z)=
((I+Z_N)^{-1}Z_1,\ldots ,(I+Z_N)^{-1}Z_{N-1},
i(I+Z_N)^{-1}(I-Z_N)),$$
for $Z=\left(Z_1\ldots Z_N\right)\in \cB _N(\cE )$. Note that since
$\left(Z_1\ldots Z_N\right)\in \cB _N(\cE )$,
each $Z_k$, $1\leq k\leq N$, is a strict contraction ($\|Z_k\|<1$), hence
$\cC (Z)$ is well-defined.

\begin{proposition}\label{Cayley}
$\cC $ is a bijection from $\cB _N(\cE )$
onto $\cG _N(\cE )$.
\end{proposition}
\begin{proof}
For $Z\in \cB _N(\cE )$ we define 
$W_k=(I+Z_N)^{-1}Z_k$ for $1\leq k\leq N-1$, and 
$W_N=i(I+Z_N)^{-1}(I-Z_N))$. Then
$$
\begin{array}{rcl}
\frac{1}{2i}(W_N-W^*_N)&=&(I+Z_N)^{-1}(I-Z_NZ^*_N)(I+Z^*_N)^{-1}\\
 & & \\
 &>&(I+Z_N)^{-1}(Z_1Z^*_1+\ldots +Z_{N-1}Z^*_{N-1})(I+Z^*_N)^{-1}\\
 & & \\
 &=&W_1W^*_1+\ldots +W_{N-1}W^*_{N-1},
\end{array}
$$
so that $\cC (Z)\in  \cG _N(\cE )$.

If $W=\left(W_1 \ldots  W_N\right) \in \cG _N(\cE )$, then 
the imaginary part of $W_N$ is strictly positive, therefore
$i+W_N$ is invertible and $\cC $ is one-to-one. Also, one easily
verifies that for $W\in \cG _N(\cE )$,
$$\cC ^{-1}(W)=
\left(2i(i+W_N)^{-1}W_1,\ldots ,2i(i+W_N)^{-1}W_{N-1},
(i+W_N)^{-1}(i-W_N)\right).
$$
\end{proof}

We now introduce the Szeg\"o kernel of $\cG _N(\cE )$
by the formula:
$$K_{\cG }(W,W')=F(W)F(W')^*, \quad W,W'\in \cG _N(\cE ),$$
where
$$F(W)=2E(\cC ^{-1}(W))\left((i+W_N)^{-1}\right)^{\oplus \infty }.
$$

\begin{lemma}\label{luni}
$(a)$  $K_{\cG }$ is a positive definite kernel on  
$\cG _N(\cE )$.

\smallskip
\noindent
$(b)$ The set $\{F(W)^*\cE \mid W\in \cG _N(\cE )\}$
is total in $\oplus _{k\geq 0}\cE _k$.

\smallskip
\noindent
$(c)$ For any $T\in \cL (\cE )$ and $W,W'\in \cB _N(\cE )$, 
$$F(W)(\frac{1}{2i}(W_NT-T(W')^*_N)-
\sum _{k=1}^{N-1}W_kT(W')^*_k)^{\oplus \infty }F(W')^*=T.$$
\end{lemma}
\begin{proof}
$(a)$ and $(b)$ follow directly from Lemma ~\ref{duminica}.
For $(c)$ we notice that if $W=\cC (Z)$ and $W'=\cC (Z')$
for $Z,Z'\in \cB _N(\cE )$, then 
$$T-\sum _{k=1}^NZ_kT(Z')^*_k=
4(i+W_N)^{-1}(\frac{1}{2i}(W_NT-T(W')^*_N)-
\sum _{k=1}^{N-1}W_kT(W')^*_k)(-i+(W')^*_N)^{-1},
$$
and together with Lemma ~\ref{duminica}$(c)$, this concludes
the proof. 
\end{proof}

We can define 
$$\cS =\{s_f:\cG _N(\cE )\rightarrow \cE \mid 
s_f(W)=F(W)f, \,\, f\in \oplus _{k\geq 0}\cE _k\},$$
and with the inner product 
$$\langle s_f,s_h\rangle _{\cS }=
\langle f,h\rangle _{\oplus _{k\geq 0}\cE _k},$$
$\cS$ becomes  Hilbert space. 
The connection between the spaces $\cR $ and $\cS $
is consistent with the connection between the 
Hardy spaces $H^2$ and $H^2(\Im z>0)$ on the unit disk
and, respectively, the upper half plane, to which 
$\cR $ and $\cS $ reduce for $N=1$ and $\cE =\CC $. 

We now let $P=\sum _{\sigma \in \FF _N^+}c_{\sigma }Y_{\sigma }
\in \cP ^0_N$ take values on $\cG _N(\cE )$
by the the formula:
$$
P(W)=\sum _{\sigma \in \FF _N^+}c_{\sigma }W_{\sigma },
\quad 
W\in \cG _N(\cE ).
$$

Let $\phi $ be a strictly positive functional on $\cO ^0_N$ 
and let $\{\varphi _{\sigma }\}_{\sigma \in \FF _N^+}$
be the associated orthogonal polynomials. Then 
$$K_n(W,W')=\sum _{|\sigma |\leq n}
\varphi _{\sigma }(W)\varphi _{\sigma }(W')^*,$$
is a positive definite kernel on $\cG _N(\cE )$, called
the {\em Christoffel-Darboux kernel}, 
and the Christoffel-Darboux formula is supposed to 
provide a connection
between $K_n$ and $K_{\cG }$.

One of the reasons for the interest in $K_n$ in the classical, 
one variable case, is that $K_n$ is a reproducing kernel 
for the set of polynomials of degree at most $n$ with respect
to the inner product induced by $\phi $. A similar result can be 
obtained in the non-commutative case.
Thus, we have to consider $\cP ^0_N(\cE)$, the set 
of elements of the form $P=\sum _{\sigma \in \FF _N^+}c_{\sigma }
Y_{\sigma }$, with only finitely many $c_{\sigma }\in \cL (\cE )$ 
different from the zero operator on $\cE $. 
On  
$\cP ^0_N(\cE)$ we introduce the 
product $\sum _{\sigma \in \FF _N^+}c_{\sigma }
Y_{\sigma }\sum _{\tau \in \FF _N^+}d_{\tau }
Y_{\tau }=\sum _{\sigma ,\tau \in \FF _N^+}c_{\sigma }
d_{\tau }Y_{\sigma \tau}$ and extend $\phi $
to $\cP ^0_N(\cE)$ by setting
$\phi (\sum _{\sigma \in \FF _N^+}c_{\sigma }
Y_{\sigma })=\sum _{\sigma \in \FF _N^+}c_{\sigma }
\phi (Y_{\sigma })$. Also, extend 
$\cI $ by $\cI (\sum _{\sigma \in \FF _N^+}c_{\sigma }
Y_{\sigma })=\sum _{\sigma \in \FF _N^+}c^*_{\sigma }
\cI (Y_{\sigma })$ and define
$$\langle P,Q\rangle _{\phi }=\phi (\cI (Q)P)$$
for $P,Q\in \cP ^0_N(\cE)$. 
Let $P=\sum _{|\sigma |\leq n}c_{\sigma }
Y_{\sigma }\in \cP ^0_N$, then $P=\sum _{|\sigma |\leq n}
a_{\sigma }\varphi _{\sigma }$, where 
$a_{\sigma }=\langle P,\varphi _{\sigma }\rangle _{\phi }$, 
so that 
$K_n$ has the following reproducing property:
$$\begin{array}{rcl}
P(W')&=&\sum _{|\sigma |\leq n}
a_{\sigma }\varphi _{\sigma }(W') \\
 & & \\
 &=&\sum _{|\sigma |\leq n}\langle P,\varphi _{\sigma }\rangle _{\phi }
\varphi _{\sigma }(W') \\
 & & \\
 &=&\langle P,\sum _{|\sigma |\leq n}\varphi _{\sigma }
\varphi _{\sigma }(W')^*\rangle _{\phi } 
\\
 & & \\
 &=&
\langle P,K_{n,W'}\rangle _{\phi },
\end{array}
$$
where $K_{n,W'}(W)=K_n(W,W')$ and $W,W'\in \cG _N(\cE )$.

We can obtain now a Cristoffel-Darboux type formula.
It is not as simple and relevant as in the commutative case.

\begin{theorem}\label{CD}
The Christoffel-Darboux kernel and the Szeg\"o kernel
are related by the formula:
$$\begin{array}{rcl}
K_n(W,W')&=&F(W)\left(\frac{1}{2i}\left(\Phi _{n+1}(W)B_{n,N}\Phi _n(W')^*-
\Phi _{n}(W)B^*_{n,N}\Phi _{n+1}(W')^*\right)\right)F(W')^* \\
 & & \\
 & & -\sum _{k=1}^{N-1}F(W)W_kK_n(W,W')(W')^*_kF(W')^*.
\end{array}
$$ 
\end{theorem}
\begin{proof}
Using the recurrence equations for the orthogonal polynomials
we deduce
$$\begin{array}{rcl}
 W_NK_n(W,W') &-& K_n(W,W')(W'_N)^* \\
 & & \\
 &=&\sum _{k=0}^n\Phi _{k+1}(W)B_{k,N}\Phi _{k}(W')^*+
\sum _{k=0}^n\Phi _{k-1}(W)B^*_{k-1,N}\Phi _{k}(W')^* \\
 & & \\
 & & -\sum _{k=0}^n\Phi _{k}(W)B^*_{k,N}\Phi _{k+1}(W')^*-
\sum _{k=0}^n\Phi _{k}(W)B_{k-1,N}\Phi _{k-1}(W')^* \\
 & & \\
 &=&\Phi _{n+1}(W)B_{n,N}\Phi _{n}(W')^*-
\Phi _{n}(W)B^*_{n,N}\Phi _{n+1}(W')^*.
\end{array}
$$
The proof is concluded by an application of Lemma ~\ref{luni} $(c)$.
\end{proof}

\section{Jacobi matrices}
Let $\phi $ be a strictly positive functional on $\cO ^0_N(\cA _2)$.
In addition to the Hilbert space $\cH _{\phi }$, the GNS 
construction produces a representation of $\cO ^0_N(\cA _2)$
by operators on $\cH _{\phi }$. In this section we analyse 
this representation in some details by showing the connection 
with certain matrices of Jacobi type. Let 
$\{\varphi _{\sigma }\}_{\sigma \in \FF ^+_N}$ be the family
of orthogonal 
polynomials associated to $\phi $. Then we define
for $P\in \cO ^0_N(\cA _2)(=\cP ^0_N)$,
\begin{equation}\label{rep}
\Psi _{\phi }(P)\varphi _{\sigma }=P\varphi _{\sigma }.
\end{equation}
Formula \eqref{simetrie}
shows that each $\Psi _{\phi }(P)$ is a symmetric operator
on $\cH _{\phi }$ with dense domain $\cD$, the linear space
generated by the polynomials 
$\varphi _{\sigma }$, $\sigma \in \FF ^+_N$. Also, 
for $P,Q \in \cP ^0_N$, 
$$\Psi _{\phi }(PQ)=\Psi _{\phi }(P)\Psi _{\phi }(Q),$$
and $\Psi _{\phi }(P)\cD \subset \cD$ for any $P\in \cP ^0_N$,  
hence $\Psi _{\phi }$ is an unbounded representation of
$\cO ^0_N(\cA _2)$. Also, $\phi (P)=
\langle \Psi _{\phi }(P)1,1\rangle _{\phi }$ for $P\in \cP ^0_N$.
Of special interest are the operators $\Psi _k=\Psi _{\phi }(Y_k)$,
$k=1,\ldots ,N$, since
$\Psi _{\phi }(\sum _{\sigma \in \FF ^+_N}c_{\sigma }Y_{\sigma })
=\sum _{\sigma \in \FF ^+_N}c_{\sigma }\Psi _{\phi ,\sigma }$, 
where $\Psi _{\phi ,\sigma }=\Psi _{i_1}\ldots \Psi _{i_k}$
for $\sigma =i_1\ldots i_k$.
Since each $\Psi _k $ commutes with the complex conjugation, 
it follows from von Neumann's theorem (Theorem ~X.7 in 
\cite{RS}) that each $\Psi _k $ admits self-adjoint extensions.

Let $\{e_1,\ldots ,e_N\}$ be the standard basis
of $\CC ^N$, then 
$\{e_{i_1\ldots i_k}=e_{i_1}\otimes \ldots \otimes e_{i_k}
\mid 1\leq i_1,\ldots ,i_k\leq N\}$ is an orthonormal basis 
of the Fock space $\cF (\CC ^N)$. Let $W$ be the unitary 
operator
from $\cF (\CC ^N)$ onto $\cH _{\phi }$ such that
$W(e_{\sigma })=\varphi _{\sigma }$, $\sigma \in \FF ^+_N$.
We see that $W^{-1}\cD $ is the linear space $\cD _0$
generated by $e_{\sigma }$, $\sigma \in \FF ^+_N$, 
so that we can define
$$J_k=W^{-1}\Psi _{\phi ,k}W,\quad k=1,\ldots ,N,$$ 
on $\cD _0$. Each $J_k$ is a symmetric operator on $\cD _0$
and by Theorem ~\ref{T3}, the matrix of $J_k$ with respect
to the orthonormal basis 
$\{e_{\sigma }\}_{\sigma \in \FF ^+_N}$ is
$$J_k=\left[
\begin{array}{cccc}
A_{0,k} & B^*_{0,k} & 0 & \ldots \\
 & & & \\
B_{0,k} & A_{1,k} & B^*_{1,k} & \\
 & & & \\ 
0 & B_{1,k} & A_{2,k} & \ddots \\
 & & & \\ 
\vdots & & \ddots & \ddots 
\end{array}
\right].$$
We call $(J_1\ldots J_N)$ a Jacobi $N$-family on $\cD _0$. We can state now
the main result of this section about the modeling
of Jacobi $N$-families.

\begin{theorem}\label{Jacobi}
Let $(J_1\ldots J_N)$ a Jacobi $N$-family on $\cD _0$ such that 
$B_n=\left[B_{n,1} \ldots B_{n,N}\right]$ 
is an upper triangular invertible matrix.
Then there exists a unique strictly positive
functional $\phi $ on $\cO ^0_N(\cA _2)$ such that the map
$$W(e_{\sigma })=\varphi _{\sigma },\quad \sigma \in \FF ^+_N,$$
extends to a unitary operator and 
$$J_k=W^{-1}\Psi _{\phi ,k}W, \quad k=1,\ldots ,N.$$
\end{theorem}
\begin{proof}
Since $B_n$, $n\geq 0$, are invertible matrices, we can uniquely
determine the elements $\varphi _{\sigma }$, $\sigma \in \FF ^+_N$,
in $\cP ^0_N$ such that for $k=1,\ldots ,N$, 
$$
Y_k\varphi _{\emptyset }=
\left[\varphi _{\sigma }\right]_{|\sigma |=1}B_{0,k}+
\varphi _{\emptyset }A_{0,k},
$$
and for 
$k=1,\ldots ,N$, $n\geq 1$,
$$
Y_k\left[\varphi _{\sigma }\right]_{|\sigma |=n}=
\left[\varphi _{\sigma }\right]^T_{|\sigma |=n+1}B_{n,k}+
\left[\varphi _{\sigma }\right]^T_{|\sigma |=n}A_{n,k}+
\left[\varphi _{\sigma }\right]^T_{|\sigma |=n-1}B^*_{n-1,k}.
$$
Theorem ~\ref{T4} gives a unique strictly positive functional
$\phi $  on $\cO ^0_N(\cA _2)$ such that 
$\{\varphi _{\sigma }\}_{\sigma \in \FF ^+_N}$
is the family of orthogonal polynomials associated to $\phi $.
The GNS construction for this $\phi $ will produce the required
$W$ and $\Psi _{\phi ,k}$, as explained above.
\end{proof}

We conclude with an application concerning the moments
of positive functionals. 
Let $\phi $ be a positive functional on $\cO ^0_N(\cA _2)$.
The numbers
\begin{equation}\label{moments}
s_{\sigma }=\phi (Y_{\sigma }),\quad  \sigma \in \FF ^+_N,
\end{equation}
are called the {\em moments } of $\phi $. 
A {\em Hamburger type problem} would be to determine
conditions on the family $\{s_{\sigma }\}_{\sigma \in \FF ^+_N}$
so that there exists a positive functional on 
$\cO ^0_N(\cA _2)$ satisfying \eqref{moments}.
A solution to this problem can be obtained as follows.

\begin{theorem}\label{Ham}
The numbers $s_{\sigma }$, $\sigma \in \FF ^+_N$,
are the moments of a positive functional on 
$\cO ^0_N(\cA _2)$ if and only if 
$K(\sigma ,\tau )=s_{I(\sigma )\tau}$, 
$\sigma ,\tau \in \FF ^+_N$,
is a positive definite kernel on $\FF ^+_N$.
\end{theorem}
\begin{proof}
The classical approach extends to this setting. Thus, we notice first
that 
\begin{equation}\label{ultima}
\begin{array}{rcl}
K(\alpha \sigma ,\tau )&=&s_{I(\alpha \sigma )\tau } \\
 & & \\
 &=&s_{I(\sigma )I(\alpha )\tau } \\
 & & \\
 &=&K(\sigma I(\alpha ),\tau ),
\end{array}
\end{equation}
for $\alpha ,\sigma ,\tau \in \FF ^+_N$.
Define the sesquiliniar form on $\cP^0_N$
by the formula
$$\langle \sum _{\sigma \in \FF ^+_N}c_{\sigma }Y_{\sigma },
\sum _{\tau \in \FF ^+_N}d_{\tau }Y_{\tau }\rangle
=\sum _{\sigma ,\tau \in \FF ^+_N}
s_{I(\sigma )\tau }c_{\tau }\overline{d}_{\sigma }.
$$
Let $\cN =\{P\in \cP^0_N\mid \langle P,P\rangle =0\}$
and let $\cH $ be the Hilbert space obtained by completing
$\cP^0_N/\cN $ in the inner product 
$\langle \cdot ,\cdot \rangle $. Define the
maps $\Psi _k:\cP^0_N\rightarrow \cP^0_N$, $k=1,\ldots ,N$,
by
$$\Psi _k(\sum _{\sigma \in \FF ^+_N}c_{\sigma }Y_{\sigma })=
\sum _{\sigma \in \FF ^+_N}c_{\sigma }Y_{k\sigma }.
$$
From \eqref{ultima}
we deduce that each $\Psi _k$, $k=1,\ldots ,N$, is a 
symmetric operator
and $\Psi _k\cN \subset \cN $.
Thus, each $\Psi _4$ lifts to a symmetric
operator, still denoted $\Psi _k$, with domain 
$\cP^0_N/\cN $ in $\cH $.
Also $\Psi _k(\cP^0_N/\cN)\subset \cP^0_N/\cN $, 
so that $\Psi _{\sigma }=\Psi _{i_1}\ldots \Psi _{i_k}$
is defined for any $\sigma =i_1\ldots i_k\in \FF ^+_N$ and we can 
take
$$\phi (\sum _{\sigma \in \FF ^+_N}c_{\sigma }Y_{\sigma })=
\langle \sum _{\sigma \in \FF ^+_N}c_{\sigma }\Psi _{\sigma }
\hat{1},\hat{1}\rangle ,$$
where $\hat{1}$ denotes the class of $1\in \cP^0_N$ in 
$\cP^0_N/\cN $.
One can easily check that $\phi $ is a solution of the Hamburger 
problem with data $\{s_{\sigma }\}_{\sigma \in \FF ^+_N}$.
\end{proof}

A more specialized version of this result was proved in \cite{Ha}.

\end{document}